\documentclass{amsart}
\input{psfig.sty}
\input{epsf.tex}
\usepackage{amsfonts,amssymb,verbatim,amsmath,amsthm,latexsym,textcomp,amscd}
\usepackage{latexsym,amsfonts,amssymb,epsfig,verbatim}
\usepackage{amsmath,amsthm,amssymb,latexsym,graphics,textcomp}
\usepackage{graphicx}
\usepackage{color}
\usepackage{url}
\usepackage{epsfig}
\input xy
\xyoption{all}

\begin{document}

\newtheorem{theorem}{Theorem}[section]
\newtheorem{prop}[theorem]{Proposition}
\newtheorem{lemma}[theorem]{Lemma}
\newtheorem{cor}[theorem]{Corollary}
\newtheorem{definition}[theorem]{Definition}
\newtheorem{conj}[theorem]{Conjecture}
\newtheorem{rmk}[theorem]{Remark}
\newtheorem{claim}[theorem]{Claim}
\newtheorem{defth}[theorem]{Definition-Theorem}

\newcommand{\C}{{\mathbb C}}
\newcommand{\integers}{{\mathbb Z}}
\newcommand{\natls}{{\mathbb N}}
\newcommand{\ratls}{{\mathbb Q}}
\newcommand{\R}{{\mathbb R}}
\newcommand{\proj}{{\mathbb P}}
\newcommand{\lhp}{{\mathbb L}}
\newcommand{\tube}{{\mathbb T}}
\newcommand{\cusp}{{\mathbb P}}
\newcommand\AAA{{\mathcal A}}
\newcommand\BB{{\mathcal B}}
\newcommand\CC{{\mathcal C}}
\newcommand\DD{{\mathcal D}}
\newcommand\EE{{\mathcal E}}
\newcommand\FF{{\mathcal F}}
\newcommand\GG{{\mathcal G}}
\newcommand\HH{{\mathcal H}}
\newcommand\II{{\mathcal I}}
\newcommand\JJ{{\mathcal J}}
\newcommand\KK{{\mathcal K}}
\newcommand\LL{{\mathcal L}}
\newcommand\MM{{\mathcal M}}
\newcommand\NN{{\mathcal N}}
\newcommand\OO{{\mathcal O}}
\newcommand\PP{{\mathcal P}}
\newcommand\QQ{{\mathcal Q}}
\newcommand\RR{{\mathcal R}}
\newcommand\SSS{{\mathcal S}}
\newcommand\TT{{\mathcal T}}
\newcommand\UU{{\mathcal U}}
\newcommand\VV{{\mathcal V}}
\newcommand\WW{{\mathcal W}}
\newcommand\XX{{\mathcal X}}
\newcommand\YY{{\mathcal Y}}
\newcommand\ZZ{{\mathcal Z}}

\title{Algebraic de Rham cohomology of log-Riemann surfaces of finite type}

\author{Kingshook Biswas}
\address{RKM Vivekananda University, Belur Math, WB-711 202, India}

\begin{abstract}
Log-Riemann surfaces of finite type are certain branched coverings, in a generalised sense, of $\C$.
They have finitely generated fundamental group and finitely many ramification points, but ramification points
of infinite order are allowed (in which case the map has infinite degree). Biswas and Perez-Marco
showed that any such log-Riemann surface is given by a meromorphic function $\pi$ on a punctured Riemann surface $S'$ (i.e.
a compact Riemann surface $S$ minus finitely many punctures) such that the differential $d\pi$ has
exponential singularities at the punctures (for example in the case of genus zero with
one puncture these are given by functions of the form $\pi = \int R(z) e^{P(z)} dz$, where
$R$ is a rational function, and $P$ is a polynomial of degree equal to the number of
infinite order ramification points). Let $\RR$ be the set of infinite order ramification
points and let $S^* = S' \sqcup \RR$. We define an algebraic de Rham cohomology group
$H^1_{dR}(S^*)$ (given by certain differentials with exponential singularities modulo
differentials of certain functions with exponential singularities) for which we show that
integration along curves gives a nondegenerate pairing with the relative homology group $H_1(S^*, \RR ; \C)$.
In particular the dimension of $H^1_{dR}(S^*)$ is $2g + \#\RR + (n-2)$, where $g$ is the genus
of $S$ and $n$ the number of punctures. The periods of these differentials along curves
joining the infinite order ramification
points can be written in the form $\int_{\gamma} R_1(z,w) e^{\int R_2(z,w) dz} dz$, where $z,w$ are
meromorphic functions generating the function field of $S$, $R_1, R_2$ are rational functions of $z,w$,
and $\gamma$ is a curve such that $z \to \infty$ along certain directions at the two endpoints of $\gamma$.
\end{abstract}

\bigskip

\maketitle

\tableofcontents

\section{Introduction}

\medskip

We recall that for a nonsingular projective algebraic curve $S$ over $\C$ (Riemann surface), the algebraic de Rham cohomology group
$H^1_{dR}(S)$ may be defined as the quotient of the space of differentials of the second kind (meromorphic $1$-forms
with all residues equal to zero) by the subspace of exact differentials (differentials of meromorphic functions).
Integration along closed curves gives a
nondegenerate pairing of $H^1_{dR}(S)$ with the homology $H_1(S, \C)$ of $S$,
so that the dimension of $H^1_{dR}(S)$ is $2g$, where $g$ is the genus of $S$.

\medskip

A choice of nonconstant meromorphic function $z$ on $S$ realizes $S$ as a finite sheeted branched
covering of the Riemann sphere $\hat{\C}$. In this article we define an algebraic
de Rham cohomology group and a period pairing for certain branched coverings, in a generalized sense, of $\C$, namely
{\it log-Riemann surfaces of finite type}, which are given by transcendental functions of
infinite degree. A log-Riemann surface consists of a Riemann surface together with a local holomorphic
diffeomorphism $\pi$ from the surface to $\C$ such that the set of points $\RR$ added to the surface,
when completing it with respect to the path-metric
induced by the flat metric $d\pi$, is discrete. Log-Riemann surfaces were
defined and studied in \cite{bpm1} (see also \cite{bpm4}), where it was shown that
the map $\pi$ restricted to any small enough punctured
metric neighbourhood of a point $w^*$ in $\RR$ gives a covering of a punctured disc in $\C$,
and is thus equivalent to
either $(z \mapsto z^n)$ restricted to a punctured disc $\{ 0 < |z| < \epsilon \}$
(in which case we say $w^*$ is a ramification point of order $n$) or to $(z \mapsto e^z)$
restricted to a half-plane $\{ \Re z < C \}$
(in which case we say $w^*$ is a ramification point of infinite order).

\medskip

A log-Riemann surface is said to be of finite type if it has finitely many ramification points
and finitely generated fundamental group. We will only consider those for which the set of
infinite order ramification points is nonempty (otherwise the map $\pi$ has finite degree
and is given by a meromorphic function on a compact Riemann surface). In \cite{bpm2}, \cite{bpm3},
uniformization theorems were proved for log-Riemann surfaces of finite type, which imply that a
log-Riemann surface of finite type is given by a pair $(S' = S - \{p_1, \dots, p_n \}, \pi)$, where
$S$ is a compact Riemann surface, and $\pi$ is a meromorphic function on the punctured surface $S'$
such that the differential $d\pi$ has essential
singularities at the punctures of a specific type, namely {\it exponential singularities}.

\medskip

Given a germ of meromorphic function $h$ at a point $p$ of a Riemann surface, a function $f$ with an isolated singularity
at $p$ is said to have an exponential singularity of type $h$ at $p$ if locally $f = g e^h$ for
some germ of meromorphic function $g$ at $p$, while a 1-form $\omega$ is said to have an
exponential singularity of type $h$ at $p$ if locally $\omega = \alpha e^h$ for
some germ of meromorphic 1-form $\alpha$ at $p$. Note that the spaces of germs of functions and 1-forms
with exponential singularity of type $h$ at $p$ only depend on the equivalence class $[h]$ in the space
$\mathcal{M}_p/\OO_p$ of germs of meromorphic functions at $p$ modulo germs of holomorphic functions at $p$.

\medskip

Thus the uniformization theorems of \cite{bpm2}, \cite{bpm3} give us
$n$ germs of meromorphic functions $h_1, \dots, h_n$ at the punctures $p_1,\dots,p_n$, with poles of orders
$d_1, \dots, d_n \geq 1$ say, such that near a puncture $p_j$ the map $\pi$ is of the form $\int g_j e^{h_j} dz$,
where $g_j$ is a germ of meromorphic function near $p_j$ and $z$ a local coordinate near $p_j$. The punctures correspond
to ends of the log-Riemann surface, where at each puncture $p_j$, $d_j$ infinite order ramification points are added
in the metric completion, so that the total number of infinite order ramification points is $\sum_j d_j$. The $d_j$
infinite order ramification points added at a puncture $p_j$ correspond to the $d_j$ directions of approach to the
puncture along which $\Re h_j \to -\infty$ so that $e^{h_j}$ decays exponentially and $\int^z g_j e^{h_j} dz$ converges.
In the case of genus zero with one
puncture for example, which is considered in \cite{bpm2}, $\pi$ must have the form $\int R(z) e^{P(z)} dz$
where $R$ is a rational function and $P$ is a polynomial of degree equal to the number of infinite
order ramification points.

\medskip

For a log-Riemann surface of finite type we would like to define a
 period pairing between differentials and curves on the surface which takes into account the extra
 structure coming from the infinite order ramification points. It is natural to consider then,
 in addition to closed curves, also curves joining the
infinite order ramification points. As the infinite order ramification points correspond to
directions of approach to the punctures, we need to consider differentials whose integrals along
curves asymptotic to these directions converge.

\medskip

Let $S^* := S' \sqcup \RR$ be the space obtained by adding the infinite order
ramification points to the punctured surface $S'$. Let $\Omega(S^*)$ be the space of $1$-forms meromorphic
on $S'$ and with exponential singularities at the
punctures $p_1,\dots,p_n$ of types $h_1,\dots,h_n$, and let $\Omega_{II}(S^*)$ be the subspace
of such $1$-forms whose residues at all points of $S'$ vanish (this will be the analogue
of the differentials of the second kind). Let $\MM(S^*)$ be the space of functions $f$
which are meromorphic on $S'$ with exponential singularities at the
punctures $p_1,\dots,p_n$ of types $h_1,\dots,h_n$. For $f$ in $\MM(S^*)$ the differential $df$
lies in $\Omega_{II}(S^*)$, and $f$ vanishes at all the infinite
order ramification points, thus so do the integrals of $df$ over curves joining the
infinite order ramification points. We let $d\MM(S^*) \subset \Omega_{II}(S^*)$ be the space of such
differentials, and define the algebraic de Rham cohomology group of the log-Riemann surface to be
the quotient space
$$
H^1_{dR}(S^*) := \Omega_{II}(S^*) / d\MM(S^*)
$$

The integrals of cohomology classes along closed curves in $S$ and along curves joining infinite order
ramification points are then well-defined, giving a bilinear pairing between $H^1_{dR}(S^*)$ and the relative
homology $H_1(S^*, \RR ; \C)$.

\medskip

\begin{theorem} \label{mainthm1} The period pairing $H^1_{dR}(S^*) \times H_1(S^*, \RR ; \C) \to \C$
is nondegenerate. In particular
$$
dim H^1_{dR}(S^*) =  dim H_1(S^* , \RR ; \C) = 2g + (n-1) + \sum_j (d_j-1) + (n-1) = 2g + \#\RR + (n-2)
$$
\end{theorem}

\medskip

We can also consider the subspaces $\OO(S^*) \subset \MM(S^*)$ and $\Omega^0(S^*) \subset \Omega_{II}(S^*)$
of functions and $1$-forms respectively in these spaces
which are holomorphic on $S'$. The differential of any $f$ in $\OO(S^*)$ lies in $\Omega^0(S^*)$ and
we can define another cohomology group $H^1_{dR, 0}(S^*) := \Omega^0(S^*) / d\OO(S^*)$.

\medskip

\begin{theorem} \label{mainthm2} The period pairing $H^1_{dR,0}(S^*) \times H_1(S^*, \RR ; \C) \to \C$
is nondegenerate. In particular
$$
dim H^1_{dR,0}(S^*) =  dim H_1(S^*, \mathcal{R} ; \C) = 2g + \#\RR + (n-2)
$$
\end{theorem}

\medskip

The proofs of the above theorems will be given in section 5.

\medskip

We note that all the spaces of functions and differentials defined above
only depend on the types $h_1, \dots, h_n$ of the exponential
singularities of $d\pi$, which are unaltered if the differential $d\pi$ is multiplied by a non-zero
meromorphic function. It is natural then to consider a less rigid structure than that of a
log-Riemann surface. Namely we declare two maps $\pi_1, \pi_2$ inducing
log-Riemann surface structures of finite type on a surface $S$ to be equivalent if the function $d\pi_1/d\pi_2$
is meromorphic on $S$, and define an {\it exp-algebraic curve} to be an equivalence class of such finite
type log-Riemann surfaces. Then exp-algebraic curves
correspond precisely to the data of a compact Riemann surface $S$ together with a Mittag-Leffler distribution of
principal parts $(h_1, \dots, h_n)$ of meromorphic functions at punctures $p_1,\dots, p_n$.

\medskip

Periods of differentials with exponential singularities have also appeared in the work of
Bloch-Esnault (\cite{blochesnault}), who define homology and cohomology groups associated to an
irregular connection on a Riemann surface, and show that an associated pairing is nondegenerate. It is
not clear however how the groups defined in \cite{blochesnault} are related to the ones defined in this
article.

\medskip

Certain functions with exponential singularities, namely the $n$-point {\it Baker-Akhiezer functions}
(\cite{baker}, \cite{akhiezer}), have been used in the algebro-geometric integration of
integrable systems (see, for example, \cite{krichever3},
\cite{krichever2} and the surveys \cite{krichever1}, \cite{dubrovin},
\cite{krichnov}, \cite{dubkrichnov}). Given a divisor $D$
on $S'$, an $n$-point Baker-Akhiezer function (with respect to the data $(\{p_j\}, \{h_j\}, D)$)
is a function $f$ in the space $\MM(S^*)$ satisfying the additional properties that the divisor $(f)$ of
zeroes and poles of $f$ on $S'$ satisfies $(f) + D \geq 0$, and that $f \cdot e^{-h_j}$ is holomorphic at $p_j$ for all
$j$. For $D$ a non-special divisor of degree at least $g$,
the space of such Baker-Akhiezer functions is known to have dimension $deg D - g + 1$.

\medskip

In section 3 we show how to naturally associate to the
data ${\mathcal{H}} = \{h_j\}$ a degree zero line bundle $L_{\mathcal{H}}$ together with
a meromorphic connection $\nabla_{\mathcal{H}}$ and a single-valued horizontal section $s_{\mathcal{H}}$.
These can be used to construct a non-zero meromorphic function $f_0$ on $S'$ with the prescribed
types of exponential singularities $h_1, \dots, h_n$ at the punctures $p_1, \dots, p_n$, so that all the spaces of functions
and $1$-forms defined above associated to the data $\{h_j\}$ are non-trivial.

\medskip

Finally we remark that functions and differentials with exponential singularities
on compact Riemann surfaces have also been studied by Cutillas
(\cite{cutillas1}, \cite{cutillas2}, \cite{cutillas3}), where they
arise naturally in the solution of the {\it Weierstrass problem} of realizing arbitrary divisors
on compact Riemann surfaces, and by Taniguchi (\cite{taniguchi1}, \cite{taniguchi2}), where entire functions satisfying certain
topological conditions (called "structural finiteness") are shown to be precisely those entire functions
whose derivatives have an exponential singularity at ${\infty}$, namely functions of the
form $\int Q(z) e^{P(z)} \ dz$, where $P, Q$ are polynomials.

\medskip

{\bf Acknowledgements.} The author would like to thank R. Perez-Marco for
numerous helpful discussions. The author was partly supported by Department of Science and
Technology research project grant DyNo. 100/IFD/8347/2008-2009.

\bigskip

\section{Log-Riemann surfaces of finite type and exp-algebraic curves}

\medskip

We recall some basic definitions and facts from \cite{bpm1}, \cite{bpm2},
\cite{bpm3}.

\medskip

\begin{definition} A log-Riemann surface is a pair $(S, \pi)$
where $S$ is a Riemann surface and $\pi : S \to \C$ is
a local holomorphic diffeomorphism such that the set of points $\RR$ added
to $S$ in the completion $S^* := S \sqcup \RR$
with respect to the path metric induced by the flat
metric $|d\pi|$ is discrete.
\end{definition}

\medskip

The map $\pi$ extends to the metric completion $S^*$ as
a $1$-Lipschitz map. In \cite{bpm1} it is shown that the map $\pi$
restricted to a sufficiently small punctured metric neighbourhood
$B(w^*, r) - \{w^*\}$ of a ramification point is a covering of
a punctured disc $B(\pi(w^*), r) - \{\pi(w^*)\}$ in $\C$, and
so has a well-defined degree $1 \leq n \leq +\infty$, called the order
of the ramification point (we assume that the order is always at least $2$,
since order one points can always be added to $S$ and $\pi$ extended
to these points in order to obtain a log-Riemann surface).

\medskip

\begin{definition} A log-Riemann surface is of finite type if it has
finitely many ramification points and finitely generated fundamental group.
\end{definition}

\medskip

For example, the log-Riemann surface given by $(\C, \pi = e^z)$ is of finite type (with the
metric $|d\pi|$ it is isometric to the Riemann surface of the logarithm, which 
has a single ramification point of infinite order), as is the log-Riemann surface
given by the Gaussian integral $(\C, \pi = \int e^{z^2} dz)$, which has two ramification points,
both of infinite order, as in the figure below:

\medskip

{\hfill {\centerline {\psfig {figure=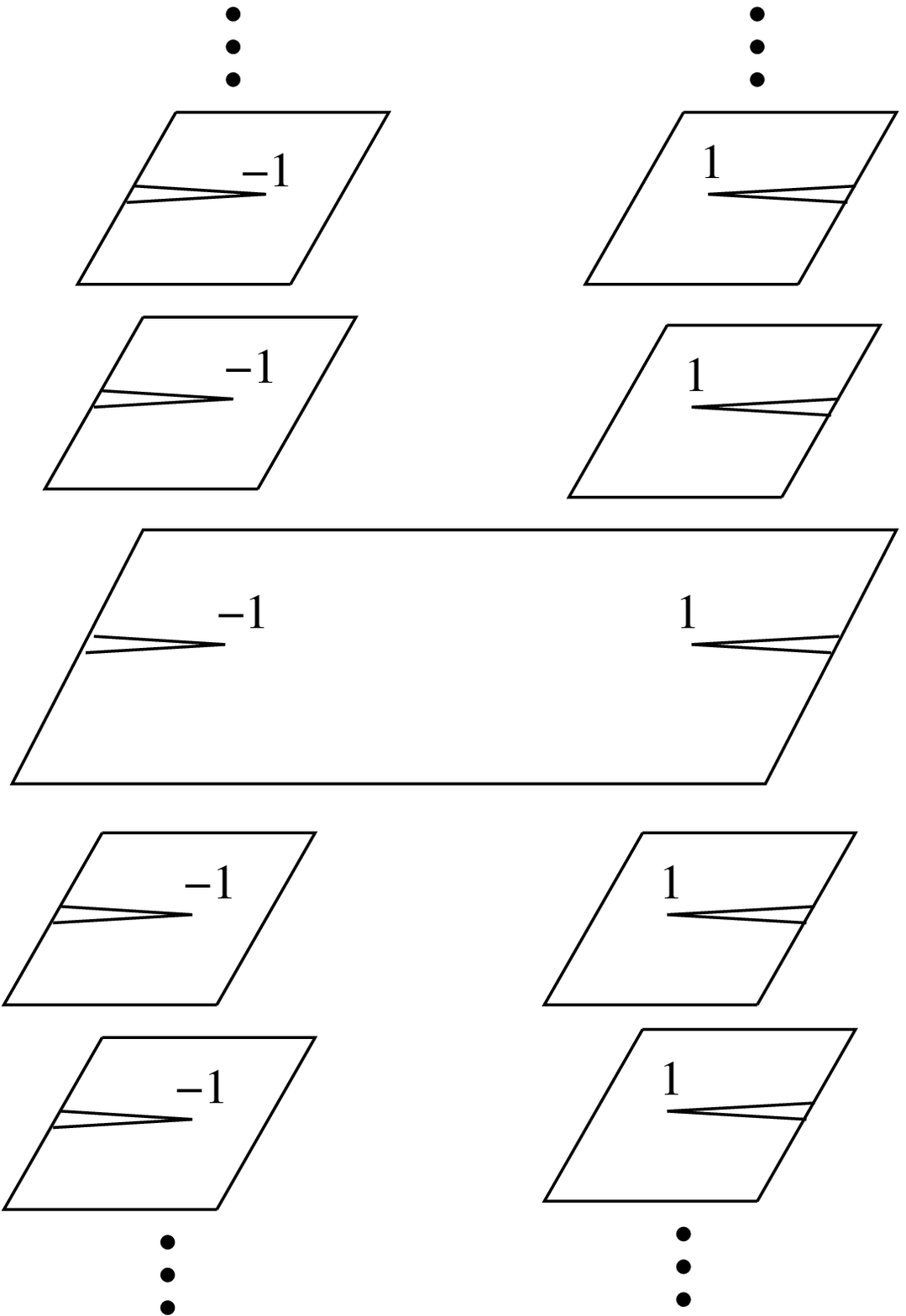,height=7cm}}}}

\medskip

\centerline{Log-Riemann surface of the Gaussian integral}

\medskip

In \cite{bpm3}, it is shown that a log-Riemann surface of finite type (which has at least one
infinite order ramification point) is of the form
$(S', \pi)$, where $S'$ is a punctured compact Riemann surface $S' = S - \{p_1, \dots, p_n\}$
and $\pi$ is meromorphic on $S'$ and $d\pi$ has exponential singularities at the punctures
$p_1, \dots, p_n$. Let $h_1, \dots, h_n$ be the types of the exponential singularities
of $d\pi$ at the punctures $p_1, \dots, p_n$.
As described in \cite{bpm3}, each puncture $p_j$ corresponds to an
end of the log-Riemann surface where $d_j$ infinite order ramification points are added,
$d_j$ being the order of the pole of $h_j$ at $p_j$.

\medskip

Let $w^*$ be an infinite order ramification point associated to a puncture $p_j$.
An $\epsilon$-ball $B_{\epsilon}$ around $w^*$ is isometric to the $\epsilon$-ball
around the infinite order ramification
point of the Riemann surface of the logarithm (given by cutting and pasting infinitely many
discs together), and there is an argument function $\arg_{w^*} : B^*_{\epsilon} \to \R$ defined
on the punctured ball $B^*_{\epsilon}$. While the function $\pi$, which is of the form
$\pi = \int e^{h_j} \alpha_j$ in a punctured neighbourhood of $p_j$ for some meromorphic $1$-form
$\alpha_j$, extends continuously to $w^*$
for the metric topology on $S^*$, in general functions of the form $f = \int e^{h_j} \alpha$
(where $\alpha$ is a $1$-form meromorphic near $p_j$) do {\bf not} extend continuously
to $w^*$ for the metric topology (\cite{bpm4}). Limits of these functions in sectors
$\{ p \in B^*_{\epsilon} : c_1 < \arg_{w^*}(p) < c_2 \}$ do exist however and are independent
of the sector; we say that the function is {\it Stolz continuous} at points of $\RR$.

\medskip

\begin{definition} Define spaces of functions and $1$-forms on $S^*$:
\begin{align*}
\MM(S^*) & := \{ f \hbox{ meromorphic function on } S' : f \hbox{ has exponential singularities at } \\
         & \quad \quad p_1, \dots, p_n \hbox{ of types } h_1, \dots, h_n \} \\
\OO(S^*) & := \{ f \in \MM(S^*) : f \hbox{ holomorphic on } S' \} \\
\Omega(S^*) & := \{ \omega \hbox{ meromorphic } 1-\hbox{form on } S' : \omega \hbox{ has exponential singularities at } \\
         & \quad \quad p_1, \dots, p_n \hbox{ of types } h_1, \dots, h_n \} \\
\Omega_{II}(S^*) & := \{ \omega \in \Omega(S^*) : \hbox{Residues of } \omega \hbox{ at all points of } S' \hbox{ vanish} \} \\
\Omega^0(S^*) & := \{ \omega \in \Omega(S^*) : \omega \hbox{ holomorphic on } S' \}
\end{align*}
\end{definition}

\medskip

Functions in $\MM(S^*)$ are Stolz continuous at points of $\RR$ taking the value $0$ there. The integrals
of $1$-forms $\omega$ in $\Omega_{II}(S^*)$ over curves $\gamma : [a, b] \to S^*$ joining points $w^*_1, w^*_2$ of
$\RR$ converge if $\gamma$ is disjoint from the poles of $\omega$ and
tends to these points through sectors $\{ p \in B^*_{\epsilon} : c_1 < \arg_{w^*_1}(p) < c_2 \},
\{ p \in B^*_{\epsilon} : c_1 < \arg_{w^*_2}(p) < c_2 \}$ (since any primitive of $\omega$ on a sector
is Stolz continuous).

\medskip

The definitions of the above spaces only depend on the types $\{ [h_i] \in \MM_{p_i}/\OO_{p_i} \}$ of the exponential
singularities of the $1$-form $d\pi$, which do not change if $d\pi$ is multiplied by a meromorphic function.
It is natural to define then a structure less rigid than that of a log-Riemann surface of finite type.

\medskip

\begin{definition}[Exp-algebraic curve] Given a punctured compact Riemann surface $S' = S - \{p_1, \dots, p_n\}$,
two meromorphic functions $\pi_1, \pi_2$ on $S'$ inducing log-Riemann surface structures of finite type
are considered equivalent if $d\pi_1/d\pi_2$ is meromorphic on the compact surface $S$. An exp-algebraic curve is an
equivalence class of such log-Riemann surface structures of finite type.
\end{definition}

\medskip

It follows from the uniformization theorem of \cite{bpm3} that
an exp-algebraic curve is given by the data of a punctured compact Riemann surface and
$n$ (equivalence classes of) germs of meromorphic functions $\mathcal{H} = \{ [h_i] \in \MM_{p_i}/\OO_{p_i} \}$
with poles at the punctures.

\medskip

We can associate
a topological space $\hat{S}$ to an exp-algebraic curve, given as a set by $\hat{S} = S' \cup \RR$,
where $\RR$ is the set of infinite ramification points added with respect to any map $\pi$ in the
equivalence class of log-Riemann surfaces of finite type, and the topology is the weakest
topology such that all maps $\tilde{\pi}$ in the equivalence class extend continuously to $\hat{S}$.

\medskip

Finally, for a meromorphic function $f$ on $S'$ (respectively meromorphic $1$-form $\omega$ on $S'$)
with exponential singularities of
types $h_1, \dots, h_n$ at points $p_1, \dots, p_n$ we can define a divisor $(f) = \sum_{p \in S} n_p \cdot p$
(respectively $(\omega) = \sum_{p \in S} m_p \cdot p$) by $n_p = ord_p(f)$ if $p \in S'$ and
$n_p = ord_{p_i}(g)$ if $p = p_i$, where $g$ is a germ of meromorphic function at $p_i$ such that $f = ge^{h_i}$
(respectively $m_p = ord_p(\omega)$ if $p \in S'$ and
$n_p = ord_{p_i}(\alpha)$ if $p = p_i$, where $\alpha$ is a germ of meromorphic $1$-form
at $p_i$ such that $\omega = \alpha e^{h_i}$).

\medskip

Note that the divisor $(f)$ can also be defined by $n_p = Res(df/f, p)$, so it follows from the Residue
Theorem applied to the meromorphic $1$-form $df/f$ that the divisor $(f)$ has degree zero.

\medskip

\section{Exp-algebraic curves and line bundles with meromorphic connections}

\medskip

Let $(S, \mathcal{H} = \{[h_1], \dots, [h_n]\})$ be an exp-algebraic curve, where $S$ is a compact Riemann surface of genus $g$
and $h_1, \dots, h_n$ are germs of meromorphic functions at points $p_1, \dots, p_n$. Let $\Omega(S)$ be the
space of holomorphic $1$-forms on $S$. The data $\mathcal{H}$ defines a degree zero line bundle $L_{\mathcal{H}}$
together with a transcendental section $s_{\mathcal H}$ of this line bundle which is non-zero on the punctured
surface $S'$ as follows:

\medskip

Solving the Mittag-Leffler problem locally for the distribution $\{h_1, \dots, h_n\}$ gives meromorphic
functions on an open cover such that the differences are holomorphic on intersections, and hence gives
an element of $H^1(S, \OO)$. Under the exponential this gives a degree zero line bundle as an element
of $H^1(S, \OO^*)$. Explicitly this is constructed as follows:

\medskip

Let $B_1, \dots, B_n$ be pairwise disjoint coordinate disks around the
punctures $p_1, \dots, p_n$ and let $V$ be an open subset of $S'$ intersecting each disk $B_i$ in
an annulus $U_i = V \cap B_i$ around $p_i$ such that $\{B_1, \dots, B_n, V \}$ is an open cover of $S$.
Define a line bundle $L_{\mathcal{H}}$ by taking the functions $e^{-h_i}$ to be the transition functions
for the line bundle on the intersections $U_i$. Define a holomorphic non-vanishing section
of $L_{\mathcal{H}}$ on $S'$ by:
$$
s_{\mathcal{H}} :=  \begin{matrix} 1 & \hbox{ on } V \\
                                    e^{-h_i} & \hbox{ on } B_i - \{p_i\} \\
                            \end{matrix} 
$$
Define a connection $\nabla_{\mathcal{H}}$ on $L_{\mathcal{H}}$ by declaring
that $\nabla_{\mathcal{H}}(s_{\mathcal{H}}) = 0$. Then for any holomorphic section
$s$ on $V$, $s = f s_{\mathcal{H}}$ for some holomorphic function $f$, and
$\nabla_{\mathcal{H}}(s) = df s_{\mathcal{H}}$, so $\nabla_{\mathcal{H}}$ is holomorphic
on $V$. On each disk $B_i$, letting $s_i$ be the section which is constant equal to $1$
on $B_i$ (with respect to the trivialization on $B_i$), for any holomorphic section
$s$ on $B_i$, $s = f s_i$ for some holomorphic function $f$, and $s_i = e^{h_i} s_{\mathcal{H}}$,
so
\begin{align*}
\nabla_{\mathcal{H}}(s) & = \nabla_{\mathcal{H}}(f s_i) \\
                        & = \nabla_{\mathcal{H}}(fe^{h_i}s_{\mathcal{H}}) \\
                        & = (df + f dh_i)e^{h_i} s_{\mathcal{H}} \\
                        & = (df + f dh_i) s_i \\
\end{align*}

thus the connection $1$-form of $\nabla_{\mathcal{H}}$ with respect to $s_i$ is given by
$dh_i$, so $\nabla_{\mathcal{H}}$ is meromorphic on $B_i$ with a single pole at $p_i$
of order $d_i + 1 \geq 2$.

\medskip

Let $s^*_{\HH}$ be the unique section of the dual bundle $L^*_{\mathcal{H}}$ on $S'$
such that $s_{\HH} \otimes s^*_{\HH} = 1$ on $S'$. Then for any non-zero meromorphic section
$s$ of $L_{\HH}$, the function $f := s \otimes s^*_{\mathcal{H}}$ is meromorphic on $S'$
with exponential singularities at $p_1, \dots, p_n$ of types $h_1, \dots, h_n$, and the
divisors of $s$ and $f$ coincide. Thus the
line bundle $L_{\mathcal{H}}$ has degree zero. In summary we have:

\medskip

\begin{theorem} \label{existence} For any log-Riemann surface of finite type $S^*$, the line bundle $L_{\mathcal{H}}$
has degree zero and the maps $s \mapsto s \otimes s^*_{\HH}, f \mapsto f \cdot s_{\HH}$
(respectively $\alpha \mapsto \alpha \otimes s^*_{\HH}, \omega \mapsto \omega \cdot s_{\HH}$)
are mutually inverse isomorphisms between the spaces of meromorphic sections of $L_{\HH}$ and
$\MM(S^*)$ (respectively the spaces of meromorphic $L_{\HH}$-valued $1$-forms and $\Omega(S^*)$)
preserving divisors.

\medskip

In particular the vector spaces $\MM(S^*), \OO(S^*), \Omega(S^*), \Omega_{II}(S^*), \Omega^0(S^*)$
are non-zero.
\end{theorem}

\medskip

\noindent{\bf Proof:} Since the isomorphisms above preserve divisors,
the spaces $\OO(S^*), \Omega^0(S^*)$ correspond to the spaces of meromorphic sections of $L_{\HH}$ and
meromorphic $L_{\HH}$-valued $1$-forms which are holomorphic on $S'$, both of which are non-empty. $\diamond$

\medskip

\begin{prop} The correspondence $\HH \mapsto (L_{\HH}, \nabla_{\HH})$ gives a one-to-one
correspondence between exp-algebraic structures on $S$ and degree zero line bundles on $S$ with
meromorphic connections with all poles of order at least two, zero residues, and trivial monodromy.
\end{prop}

\medskip

\noindent{\bf Proof:} Since the connection $1$-form of $\nabla_{\HH}$ is given by $dh_i$ on $B_i$,
all residues of $\nabla_{\HH}$ are equal to zero, while the monodromy of $\nabla_{\HH}$ is trivial
since $s_{\HH}$ is a single-valued horizontal section.

\medskip

Conversely, given such a meromorphic connection $\nabla$ on a degree zero line bundle $L$, if $p_1, \dots, p_n$
are the poles of $\nabla$ and $\omega_1, \dots, \omega_n$ are the connection $1$-forms of $\nabla$ with
respect to trivializations near $p_1, \dots, p_n$, then each $\omega_i$ has zero residue at $p_i$ and pole
order at least two, hence there exist meromorphic germs $h_1, \dots, h_n$ near $p_1, \dots, p_n$ such that
$\omega_i = dh_i$. We obtain an exp-algebraic curve $(S, \HH(L, \nabla))$.

\medskip

It is clear for an exp-algebraic curve $(S, \HH)$ that $\HH(L_{\HH}, \nabla_{\HH}) = \HH$, so the
correspondences are inverses of each other. $\diamond$

\medskip

Finally we remark that by Serre Duality, the degree zero line bundle $L_{\mathcal{H}}$, given as an element of $H^1(S, \OO)$,
can also be described as an element of $H^0(S, \Omega)^* = \Omega(S)^*$ using residues, as the linear
functional
\begin{align*}
Res_{\mathcal{H}} : \Omega(S) & \rightarrow \C \\
                      \xi     & \mapsto  \sum_i Res(\xi \cdot h_i, p_i) \\
\end{align*}

\medskip

\section{The local period pairing}

\medskip

Let $p = p_i$ be one of the punctures of the exp-algebraic curve $(S, \mathcal{H})$,
and let $(U, w = 1/z)$ be a coordinate disk near $p$ such that $z(p) = \infty$ and $h = h_i$ is given by $h = z^d$,
where $d = d_i$ is the order of the pole of $h_i$ at $p_i$.
The completion $U^*$ of the punctured disk $U' = U - \{p\}$ with respect to the metric $|e^{z^d}| |dz|$
is given by adding $d$ infinite order ramification points $w^*_0, \dots, w^*_{d-1}$ to $U'$,
$U^* = U' \cup \{w^*_1, \dots, w^*_d\}$, such that $z \in U'$ converges to $w^*_j$ if
$z \to \infty$ along any ray $\arg z = \theta$ with $|\theta - (2j-1)\pi/d| < \pi/2d$.

\medskip

\begin{definition} Let ${\OO}_p([h]), \Omega_p([h])$ be the spaces of holomorphic
functions and $1$-forms respectively on $U'$ with exponential singularities of type
$h$ at $p$.
\end{definition}

\medskip

Let $\gamma \subset U'$ be a small circle around $p$ traversed anticlockwise, and
let $\gamma_1, \dots, \gamma_{d-1} : \R \to U'$ be curves
such that $\gamma_k(t) \to w^*_0$ as $t \to -\infty$, $\gamma_k(t) \to w^*_k$ as $t \to +\infty$
(so that $\gamma_k$ is a curve joining the infinite order ramification points $w^*_0$ and $w^*_k$).
Then the curves $\gamma, \gamma_1, \dots, \gamma_{d-1}$
form a basis for $H_1(U^*, \{w^*_0, \dots, w^*_{d-1} \}; \C)$, and integration along curves gives a
pairing of this space with the space $\Omega_p([h])/d{\OO}_p([h])$, the 'local period pairing'
at $p$. We show in the following sections that this local period pairing is nondegenerate.

\medskip

First we show:

\medskip

\begin{theorem} \label{localdim} $$ \dim \Omega_p([h])/d{\OO}_p([h]) = d $$.
\end{theorem}

\medskip

\subsection{Differentials of the form $\omega = z^k e^{z^d}
dz$.}

\medskip

Let $g(z)e^{z^d}dz$ be a differential in $\Omega_p([h])$
with $g$ meromorphic at infinity.
In the case when $g$ is a polynomial we
have:

\medskip

\begin{prop} \label{poly} \noindent (1) For $d=1$, $$ \mathbb{C}[z]e^{z} dz
= d({\C}[z]e^{z}) $$ \noindent (2) For $d \geq 2$, $$ {\C}[z]e^{z^d} dz
= d\left({\C}[z]e^{z^d}\right) \bigoplus_{j =
0}^{d-2} \C \cdot z^j e^{z^d} dz $$
\end{prop}

\medskip

\noindent {\bf Proof:} For $d = 1$, $e^z dz = d(e^z)$ and for $k \geq 0$,
$z^{k+1}e^z dz = d(z^{k+1}e^z) - (k+1)z^k e^z dz$ (integration by
parts), so the result follows by induction.

Let $d \geq 2$. For $0 \leq k \leq d-2$ there is nothing to prove,
for $k = d-1$ we have $z^{d-1}e^{z^d} dz = (1/d) d(e^{z^d})$, and
for $k \geq d-1$, $z^{k+1} = z^{d-1}z^{k-d+2}$ so $$ z^{k+1}
e^{z^d} dz = {1 \over d} d\left(z^{k-d+2} e^{z^d}\right) - {k-d+2
\over d} z^{k-d+1} e^{z^d} dz $$ (integration by parts) so the
result follows by induction. $\diamond$

\medskip

For $g$ a Laurent polynomial we have:

\medskip

\begin{prop} \label{laurent} For $d \geq 1$, $$ {\C}[z, 1/z] e^{z^d} dz =
d\left({\C}[z, 1/z]e^{z^d}\right) \bigoplus_{j = -1}^{d-2} {\C} \cdot z^j e^{z^d} dz $$
\end{prop}

\noindent {\bf Proof:} Let $1 \leq k \leq d$, consider $(1/z^k)e^{z^d} dz$.
If $k=1$ there is nothing to prove. If $1 < k \leq d$, then $$
{e^{z^d} \over z^k} dz = d\left({-e^{z^d} \over (k-1)
z^{k-1}}\right) + {d \over k-1} z^{d-k} e^{z^d} dz $$ (integration
by parts) and we are done by the previous Proposition since $d - k \geq 0$. For $k > d$ the
result follows by induction since $k - d < k$. $\diamond$

\medskip

Let $\mathcal{M}_p$ be the space of germs of meromorphic functions at $p$.
Theorem \ref{localdim} follows immediately from the following proposition:

\medskip

\begin{prop}\label{localbasis} For $d \geq 1$, $$ \mathcal{M}_p \cdot e^{z^d} dz =
d(\mathcal{M}_p \cdot e^{z^d}) \bigoplus_{j = -1}^{d-2} \C \cdot {z^j e^{z^d}} dz $$
\end{prop}

\medskip

\noindent{\bf Proof:} A straightforward induction using the proof of
Proposition \ref{laurent} above shows that for $k > 1$,
$$
{e^{z^d} \over z^k} dz = \sum_{j = 1}^{k-1} \alpha_{j,k} d\left({e^{z^d} \over
z^j}\right) + \sum_{j = -1}^{d-2} \beta_{j,k} z^j e^{z^d} dz
$$
where the constants $\alpha_{j,k}, \beta_{j,k}$ are given by
$$
\displaylines{ \alpha_{j,(j + kd)+1} = \left({-1 \over
j}\right)\left({d \over j+d}\right)\left({d \over
j+2d}\right)\dots \left({d \over j+kd}\right) \ , \ j \geq 1, k
\geq 0, \cr
\beta_{d-j, d+kj} = \left({d \over j-1}\right)\left({d
\over j-1 + d}\right)\left({d \over j-1+2d}\right)\dots\left({d
\over j-1+kd}\right) \ , \ 2 \leq j \leq d+1, k \geq 0 \cr }
$$
and $\alpha_{j,k} = \beta_{j,k} = 0$ otherwise.

\medskip

Given $\omega = g(z) e^{z^d} dz \in \Omega_p([h])$ with $g$ a
convergent Laurent series $g(z) = g_n z^n + g_{n-1} z^{n-1} +
\dots + g_0 + g_{-1} z^{-1} + \dots$, by Proposition \ref{poly}
$$
(g_n z^n + \dots g_0 + g_{-1} z^{-1}) e^{z^d} dz \in d\left({\OO}_p\right)
\bigoplus_{j = -1}^{d-2} z^j e^{z^d} dz
$$
while for the remaining terms, using the above formula for $z^{-k} e^{z^d}
dz$, we have an equality of formal power series
$$ \left(\sum_{k =
2}^{\infty} {g_{-k} \over z^k} \right) e^{z^d} dz =
d\left(\left(\sum_{j = 1}^{\infty} \left(\sum_{k = j+1}^{\infty}
\alpha_{j,k} g_{-k} \right) {1 \over z^j} \right) e^{z^d} \right)
+ \sum_{j = -1}^{d-2} \left(\sum_{k = 2}^{\infty} \beta_{j,k}
g_{-k} \right) {z^j e^{z^d}} dz
$$
It suffices to show that
the series on the right above are convergent:

\medskip

Let $C, R > 0$ such that $|g_{-k}| \leq C R^k$. Choose $k_0 - 1 >
d \cdot (2R)^d$. Then it is easy to see from the above formulae
that $$\displaylines{ |\alpha_{j,k}| \leq \left({1 \over
j}\right)\left({1 \over (2R)^d}\right)^{(k - k_0)/d} \ , \ j \geq
1, k \geq k_0 \cr |\beta_{d-j,k}| \leq \left({1 \over
(2R)^d}\right)^{(k - k_0)/d} \ , 2 \leq j \leq d+1, k \geq k_0 \cr
}$$ so $$\displaylines{ \sum_{k \geq j+1} |\alpha_{j,k} g_{-k}| <
+\infty \cr \sum_{k \geq 2} |\beta_{d-j,k} g_{-k}| < +\infty \cr
}$$ and moreover for $j > k_0$, $\sum_{k \geq j+1} |\alpha_{j,k}
g_{-k}| < C' / j$ (for some $C'>0$), so the power series with
coefficient of $z^j$ equal to $\sum_{k \geq j+1} \alpha_{j,k}
g_{-k}$ has a positive radius of convergence. The proposition
follows. $\diamond$

\medskip

\subsection{Nondegeneracy of the local period pairing}

\medskip

Define a ''local period mapping'' by
\begin{align*}
 \Phi : \Omega_p([h]) & \to {\C}^d \\
                   \omega     & \mapsto \left( \int_{\gamma}
                   \omega, \int_{\gamma_1} \omega, \dots,
                   \int_{\gamma_{d-1}} \omega \right) \\
\end{align*}

\medskip

The nondegeneracy of the local period pairing follows from the following Theorem:

\medskip

\begin{theorem} \label{localperiod} For $\omega \in \Omega_p([h]), \Phi(\omega) = 0$
if and only if $\omega \in d{\OO}_p([h])$.
\end{theorem}

\medskip

For the computations which follow, it will be convenient to take the local coordinate $w = 1/z$
for this section {\bf only} such that $h(z) = -z^d$. As before we let $\gamma$ be a small circle
going anticlockwise around $p$, and let $\gamma_1, \dots, \gamma_{d-1}$ be curves joining the
ramification point $w^*_0$ to $w^*_1, \dots, w^*_{d-1}$ respectively.

\medskip

By Proposition \ref{localbasis}, it suffices to consider the periods of the differentials
$z^j e^{-z^d} dz, j=-1, \dots, d-2$. Let $\Pi$ be the ''period
matrix'' 
$$ 
\Pi = \left ( \begin{matrix} \int_{\gamma} z^{-1}e^{-z^d} dz
& \int_{\gamma} e^{-z^d} dz & \ldots &\int_{\gamma} z^{d-2}
e^{-z^d} dz \\ \int_{\gamma_1} z^{-1}e^{-z^d} dz
&\int_{\gamma_1} e^{-z^d} dz & \ldots &\int_{\gamma_1} z^{d-2}
e^{-z^d} dz \\ \vdots &\vdots &\ddots & \vdots \\
\int_{\gamma_{d-1}} z^{-1} e^{-z^d} dz & \int_{\gamma_{d-1}}
e^{-z^d} dz & \ldots &\int_{\gamma_{d-1}} z^{d-2}e^{-z^d} dz \\ \end{matrix} \right
) $$

Then Theorem \ref{localperiod} follows from:

\medskip

\begin{prop} \label{nonsingular} $\Pi$ is nonsingular.
\end{prop}

\medskip

\noindent{\bf Proof:} The first row of $\Pi$ is $(2\pi i, 0,0,\dots,0)$. It suffices
to show then that the $(d-1)$-by-$(d-1)$ minor 
$$ 
\Pi_1 = \left (\begin{matrix} \int_{\gamma_1} e^{-z^d} dz & \ldots &\int_{\gamma_1}
z^{d-2} e^{-z^d} dz \\ \vdots &\vdots &\ddots & \vdots \\
\int_{\gamma_{d-1}} e^{-z^d} dz & \ldots &\int_{\gamma_{d-1}}
z^{d-2}e^{-z^d} dz \\ \end{matrix} \right ) 
$$ 
is nonsingular. Let $\omega_k = e^{2\pi ik/d}$ and $v_k$ be the row vector 
$$ 
v_k = (\int_{0}^{\omega_k \cdot \infty} e^{-z^d} dz, \dots,
\int_{0}^{\omega_k \cdot \infty} z^{d-2} e^{-z^d} dz)
$$ 
Then the rows of $\Pi_1$ are $v_1 - v_0, v_2 - v_0, \dots, v_{d-1} - v_0$.

\medskip

A change of variables $z = \omega_k t^{1/d}$ gives 
$$
\int_{0}^{\omega_k \cdot \infty} z^j e^{-z^d} dz = {1 \over d}
{\omega_k}^{j+1} \Gamma\left({j+1 \over d}\right) 
$$ 
for $0 \leq k \leq d-1, j \geq 0$. Let $\Pi_2$ be the matrix with rows
$v_1,\dots, v_{d-1}$, then
\begin{align*}
 &\det \Pi_2 =\left |
\begin{matrix} {1 \over d} \Gamma \left ({1\over d}\right ) \omega_1& {1
\over d} \Gamma \left ({2\over d}\right ) \omega_1^2 & \ldots & {1
\over d} \Gamma \left ({d-1 \over d}\right ) \omega_1^{d-1} \cr {1
\over d} \Gamma \left ({1\over d}\right ) \omega_2 &{1 \over d}
\Gamma \left ({2\over d}\right ) \omega_2^2 & \ldots &{1 \over d}
\Gamma \left ({d-1\over d}\right ) \omega_2^{d-1} \\ \vdots
&\vdots &\ddots & \vdots \\ {1 \over d} \Gamma \left ({1\over
d}\right ) \omega_{d-1} & {1 \over d} \Gamma \left ({2\over
d}\right ) \omega_{d-1}^2 & \ldots &{1 \over d} \Gamma \left ({d-1
\over d}\right ) \omega_{d-1}^{d-1} \\ \end{matrix} \right | \\ &= \left({1
\over d}\right)^{(d-1)^2} \Gamma \left ({1\over d}\right ) \Gamma
\left ({2\over d}\right )\ldots \Gamma \left ({d-1\over d}\right )
\left | \begin{matrix} \omega_1& \omega_1^2 & \ldots & \omega_1^{d-1}
\\ \omega_2 &\omega_2^2 & \ldots & \omega_2^{d-1} \\ \vdots
&\vdots &\ddots & \vdots \\ \omega_{d-1} & \omega_{d-1}^2 &
\ldots & \omega_{d-1}^{d-1} \\ \end{matrix} \right | \\ &= \left({1 \over
d}\right)^{(d-1)^2} \Gamma \left ({1\over d}\right ) \Gamma \left
({2\over d}\right )\ldots \Gamma \left ({d-1\over d}\right )
(\omega_1 \dots \omega_{d-1}) \Pi_{1 \leq i \neq j \leq d-1}
(\omega_i - \omega_j)
\\ & \neq 0 \\
\end{align*}
 (using the formula for the Vandermonde
determinant). It is easy to see that $(-v_0) = v_1 + v_2 + \dots +
v_{d-1}$, giving
\begin{align*}
 \det \Pi_1 & = \det (v_1 - v_0,
\dots, v_{d-1} - v_0) \\
           & = \det (v_1, v_2, \dots, v_{d-1}) + \det(-v_0, v_2, \dots, v_{d-1}) +
\det(v_1, -v_0, v_3, \dots, v_{d-1}) + \\ & \dots +
\det(v_1,\dots, v_{d-2}, -v_0) \\ & = d \det \Pi_2 \\
\end{align*}
so
$\det \Pi = 2\pi i \det \Pi_1 = 2\pi i d \det \Pi_2 \neq 0$.
$\diamond$

\bigskip

\subsection{Local periods as residues of Borel transforms}

\bigskip

In this section we show that the local periods of a differential $\omega \in \Omega_p([h])$
can also be represented as residues of the Borel transform $\BB_d f$ of order $d$ of a certain function $f = f_{\omega}$
associated with $\omega$ (where $d$ is the order of the pole of $h$ at $p$), and use it to give another
proof of the fact that $\omega$ is 'exact' ($\omega \in d\OO_p([h])$) if all its local periods vanish.

\medskip

We recall the definitions
of Borel and Laplace transforms of order $d$ (\cite{balser} $\S$5.1, $\S$5.2):

\medskip

%
%
%
%
%
%
%

\medskip

\begin{definition}[Borel transform of order $d$ in the direction $\theta$] Given $d \geq 1$ and
$f$ holomorphic in a small
sector $V = \{ 0 < |w| < \rho, |\arg w - \theta| < \alpha/2 \}$
centered around the direction $\theta \in \R$ of opening $\alpha > \pi/d$
which is bounded near the origin,
let $\gamma(\theta)$ denote
a closed curve starting out from the origin along a ray $\arg w = \theta + (\pi+\epsilon)/2d$
going up to a point $w_0$ with $|w_0| < \rho$, then going clockwise
along the circle $\{ |w| = |w_0| \}$ till a ray $\arg w = \theta - (\pi-\epsilon)/2d$, and
then going back to the origin along this ray, where $\epsilon$ is chosen small enough
so that $\gamma(\theta)$ is contained in the sector $V$.

The Borel transform of order $d$ of $f$ in the direction $\theta$ is defined by
$$
(\BB_d f)(\xi) = \frac{1}{2\pi i} \int_{\gamma(\theta)} w^d f(z) \exp((\xi/w)^d) d(w^{-d})
$$
for $\xi$ in an infinite sector $W = \{ \xi \neq 0, |\arg \xi - \theta| < \epsilon/2d \}$ centered
around the direction $\theta$.

Note that for $\xi$ in $W$ the integral converges absolutely and uniformly on compacts (the exponential
term decays rapidly along the radial parts of $\gamma(\theta)$ for $\xi$ in $W$), and so $\BB_d f$ is
holomorphic in the sector $W$. Moreover by Cauchy's Theorem the integral is independent of
the choices of $w_1, \epsilon$.
\end{definition}

\medskip

\begin{definition}[Laplace transform of order $d$ in the direction $\theta$] Given $d \geq 1$
and a function $g$ holomorphic
in an narrow infinite sector $W = \{ \xi \neq 0, |\arg \xi - \theta| < \epsilon/2d \}$ centered
around the direction $\theta$ of opening $\epsilon/d < 2\pi$, such that $g$ has exponential
growth of order at most $d$ near $\infty$ ($|g(\xi)| = O(\exp(k|\xi|^d))$
 for some constant $k >0$), the Laplace transform of order $d$
of $g$ in the direction
$\theta$ is defined by
$$
(\LL_{d}g)(w) := w^{-d} \int_{0}^{\theta \cdot \infty} g(\xi) \exp(-(\xi/w)^d) d(\xi^d)
$$
for $w$ in a small sector $V = \{ 0 < |w| < \rho, |\arg w - \theta| < \alpha/2 \}$. Note
that the integral converges absolutely and uniformly on compacts (the exponential
term decays rapidly enough along the ray $[0, \theta \cdot \infty]$ for $w$ in $V$
if $\rho$ is small enough), so $\LL_{d} g$ is holomorphic in the sector $V$.
\end{definition}

\medskip

Then the following inversion theorem holds (\cite{balser}, $\S$5.3):

\medskip

\begin{theorem} \label{borellaplace} Let $f$ be holomorphic and bounded near the origin
in a small sector $V$ centered around the direction $\theta$ of opening $\alpha > \pi/d$.
Then
$$
g = \BB_d f
$$
is holomorphic and of exponential growth at most $d$ in a narrow sector
$W$ centered around the direction $\theta$, and
$$
f = \LL_d g
$$
\end{theorem}

\medskip

Let $\omega \in \Omega_p([h])$. As before, let $w = 1/z$ be a local coordinate near $p$
such that $z(p) = \infty$ and $h(z) = z^d$, and let $w^*_0, \dots, w^*_{d-1}$ be the
infinite order ramification points added when completing a punctured neighbourhood of
$p$ with respect to the metric $|e^{z^d} dz|$.

\medskip

Note that in any small sector $V$ with vertex at $w = 0$ which contains a ray towards
the origin along which $w$ tends to a ramification point $w^*_j$, we can always write $\omega$ in the
form
$$
\omega = d\left(f \cdot w^d e^{h}\right)
$$
where $f = f_{\omega}$ is holomorphic in $V$ and given by
$$
f = \frac{1}{w^d} \cdot \frac{1}{e^{h}} \cdot \int_{w^*_j}^{w} \omega = z^d \cdot \frac{1}{e^{h}} \cdot \int_{w^*_j}^{z} \omega
$$
Then $\omega \in d\OO_p([h])$ if and only if $f$ extends to a meromorphic function in a neighbourhood of $w = 0$.

\medskip

In what follows, it will be convenient to work with $f$ as a function of the $z$-variable, in which
case the sector $V$ will be a small sector with vertex at $z = \infty$. In this variable,
 if we take the sector $V$ centered around the direction $\arg z = 0$, then
 the integral defining the Borel transform is over a curve $\gamma$ going from $z = \infty$
to some point $z_1$ (with $|z_1|$ large) along a ray $\arg z = -\pi/2d - \epsilon$, then going
clockwise along the circle $\{ |z| = |z_1| \}$ up to the ray $\arg z = \pi/2d + \epsilon$ and
then going back to $z = \infty$ along this ray. Note that the curve $\gamma$ is a curve joining
$w^*_0$ to $w^*_1$. In this sector, we take $f$ to be given by
$$
f = z^d \cdot \frac{1}{e^{z^d}} \cdot \int_{w^*_1}^{z} \omega
$$

\medskip

\begin{theorem}\label{periodsresidues} Let $\omega = e^{z^d} \phi(z) \ dz \in \Omega_p([h])$
with $\phi(z) \ dz$ holomorphic near $z = \infty$. Then:

\medskip

\noindent{(1)} The Borel transform $\BB_d f$ of order $d$ in the direction $\theta = 0$
of $f_{\omega}$ extends to a meromorphic function on $\C$ which is holomorphic outside the finite set of $d$th
roots of unity $\xi = e^{2\pi i k/d}$, and has exponential growth of order at most $d$
near $\xi = \infty$.

\medskip

\noindent{(2)} The Borel transform $\BB_d f$ has at most simple poles at the $d$th roots of unity,
and the residues of $\BB_d f$ at these points are given in terms of periods of $\omega$:
$$
Res(\BB_d f, \xi = e^{2\pi i k/d}) = \left(\frac{1}{d e^{2\pi ik(d-1)/d}}\right) \left(\frac{-1}{2\pi i}\right)\int_{w^*_{k}}^{w^*_{k+1}} \omega
$$
\end{theorem}

\medskip

\noindent{\bf Proof:} We compute:
\begin{align*}
(\BB_d f)(\xi) & = \frac{1}{2\pi i} \int_{\gamma} \frac{1}{z^d} f(z) e^{(\xi z)^d} d(z^d) \\
               & = \frac{1}{2\pi i} \int_{\gamma} \left( \int_{w^*_1}^{z} \phi(t) e^{t^d} dt \right) e^{-z^d} e^{(\xi z)^d} d(z^d) \\
               & = \frac{-1}{2\pi i} \int_{\gamma} \left( \int_{z}^{w^*_1} \phi(t) e^{t^d} dt \right) e^{(\xi^d - 1)z^d} d(z^d) \\
\end{align*}
At this point we note that the inner integral from $z$ to $w^*_1$ can be taken along the curve $\gamma$ as well, so that the integral
above becomes an integral over a triangle $\{0 \leq x \leq y \leq 1\} \subset [0,1] \times [0,1]$, taking $[0,1] \times [0,1]$ to be
the domain of $\gamma \times \gamma : (x,y) \mapsto (\gamma(x), \gamma(y)) \in U^*$. Applying Fubini's Theorem gives
\begin{align*}
(\BB_d f)(\xi) & = \frac{-1}{2\pi i} \int_{\gamma} \left( \int_{w^*_0}^{t} e^{(\xi^d - 1)z^d} d(z^d) \right) \phi(t) e^{t^d} dt \\
               & = \frac{-1}{2\pi i} \int_{\gamma} \frac{1}{\xi^d - 1} e^{(\xi^d - 1)t^d} \phi(t) e^{t^d} dt \\
               & = \left(\frac{-1}{2\pi i} \int_{\gamma} \phi(t) e^{\xi^d t^d} dt\right) \cdot \frac{1}{\xi^d - 1} \\
               & = \frac{A(\xi)}{\xi^d - 1} \\
\end{align*}
where
$$
A(\xi) = \frac{-1}{2\pi i} \int_{\gamma} \phi(t) e^{\xi^d t^d} dt
$$
Now the function $A$ is holomorphic in an infinite sector centered around the positive real axis in the $\xi$-plane.
However it can be analytically continued to all of $\C^*$ by the standard device of 'rotating' the contour $\gamma$. Namely,
if $\lambda \cdot \gamma$ denotes the curve $\gamma$ multiplied by $\lambda \in S^1$, then we can extend $A$ to all of $\C^*$
by defining
$$
A(\xi) := \frac{-1}{2\pi i} \int_{(\xi/|\xi|) \cdot \gamma} \phi(t) e^{\xi^d t^d} dt
$$
for all $\xi \in \C^*$. Moreover the isolated singularity of $A$ at $\xi = 0$ is removable since $A(\xi) \to 0$ as $\xi \to 0$,
because the integrand converges to $\phi(t) dt$ and $\phi(t) dt$ is holomorphic near $t = \infty$. Thus $A$ is entire,
and it is not hard to show from the above equation that $A$ is of exponential growth of order at most $d$ near $\xi = \infty$.
This proves (1). For (2), we have
\begin{align*}
Res(\BB_d f, \xi = e^{2\pi i k/d}) & = \frac{A(e^{2\pi i k/d})}{d \cdot (e^{2\pi ik/d})^{d-1}} \\
                                   & = \left(\frac{1}{d e^{2\pi ik(d-1)/d}}\right) \left(\frac{-1}{2\pi i}\right)\int_{w^*_{k}}^{w^*_{k+1}} \omega \\
\end{align*}
$\diamond$

\medskip

As a corollary we obtain another proof of Theorem \ref{localperiod}:

\medskip

\noindent{\bf Proof of Theorem \ref{localperiod}:} Clearly all periods vanish if $\omega \in d\OO_p([h])$.
For the converse, assume all periods of $\omega$
vanish. By Proposition \ref{localbasis}, $\omega \equiv e^{z^d} P(z) dz (mod \ d\OO_p([h]))$ for some polynomial $P$ of
degree at most $d-2$, and from the proof of Proposition \ref{laurent} and the fact that $Res(\omega, p) = 0$
it follows that $e^{z^d} P(z) dz \equiv \omega_1 = e^{z^d} \phi(z) dz (mod \ d\OO_p([h]))$ for some $1$-form $\phi(z) dz$
which is holomorphic near $z = \infty$. Letting $f = f_{\omega_1}$, by the previous Theorem $\BB_d f$ is entire (since its
residues are multiples of the periods of $\omega_1$ which vanish by hypothesis) of exponential growth at most $d$, so by the
Inversion Theorem for Borel-Laplace transforms, $f = \LL_d (\BB_d f)$ defines a holomorphic germ near $z = \infty$,
thus $\omega \equiv \omega_1 \equiv 0 (mod \ d\OO_p([h])$. $\diamond$

\medskip

\section{The global period pairing}

\bigskip

Let $(S, \mathcal{H} = \{h_j\})$ be an exp-algebraic curve, where $h_1, \dots, h_n$ are germs of meromorphic
functions with poles of order $d_1, \dots, d_n \geq 1$ at points $p_1, \dots, p_n$ and $S$ is a
compact Riemann surface of genus $g$. Let $S^* = S \sqcup \RR$ be the completion of $S' = S - \{p_1, \dots, p_n\}$
with respect to some log-Riemann surface structure of finite type given by a meromorphic map $\pi : S' \to \hat{\C}$
such that $d\pi$ has exponential singularities at $p_1, \dots, p_n$ of types $h_1, \dots, h_n$.

\medskip

For each puncture $p_i$, let $w^{*(i)}_j, 0 \leq j \leq d_i - 1$ be the infinite order ramification points
associated to the puncture $p_i$, so that $\RR = \{ w^{*(i)}_j, 0 \leq j \leq d_i - 1, 1 \leq i \leq n \}$.
Let $\gamma^{(i)}$ be a small circle around $p_i$, and let
$\{ \gamma^{(i)}_j : 1 \leq j \leq d_i - 1 \}$ be curves, as in the previous section, joining the
ramification point $w^{*(i)}_0$ to the ramification points $w^{*(i)}_j, 1 \leq j \leq d_i - 1$ and
disjoint from $\gamma^{(i)}$.

\medskip

We fix a standard homology basis $a_1, \dots, a_g, b_1, \dots, b_g$ of $S$, choosing these closed curves so that they
are disjoint from the punctures $p_1, \dots, p_n$ and intersect at only one point $q_0$, so that cutting the surface
open along these curves gives a polygon $P$ with $4g$-sides containing the points $p_1, \dots, p_n$ in its interior.

\medskip

Let $\beta_i, i = 2,\dots, n,$ be curves joining $w^{*(1)}_0$ to $w^{*(i)}_n, i = 2, \dots, n$, such that the
curves $\beta_i$ are contained in the interior of the polygon $P$
(and so are disjoint from $a_1, \dots, a_g, b_1, \dots, b_g$), any two such curves intersect only at $w^{*(1)}_0$,
the curves $\beta_i$ intersect the curves $\{ \gamma^{(i)}_j \}$ only at their endpoints $w^{*(1)}_0, w^{*(i)}_n$,
and only intersect the curves $\{ \gamma^{(i)} \}$ at one point each, in $\gamma^{(1)}$ and $\gamma^{(i)}$.

\medskip

The relative homology $H_1(S^*, \RR ; \C)$ is generated by the curves in the
union $\CC := \CC_1 \cup \CC_2 \cup \CC_3 \cup \CC_4$ of the following four finite families of curves, modulo
the single relation $\gamma^{(1)} + \dots + \gamma^{(n)} = 0$:

\medskip

\noindent{1)} The $2g$ closed curves $\CC_1 := \{a_1, \dots, a_g, b_1, \dots, b_g \}$.

\medskip

\noindent{2)} The $(n-1)$ curves $\CC_2 := \{ \beta_i : 2 \leq i \leq n \}$.

\medskip

\noindent{3)} The $\sum_i (d_i - 1)$ curves $\CC_3 := \{ \gamma^{(i)}_j : 1 \leq i \leq n, 1 \leq j \leq d_i - 1 \}$.

\medskip

\noindent{4)} The $n$ closed curves $\CC_4 := \{ \gamma^{(i)} : 1 \leq i \leq n \}$.

\medskip

These curves give a well-defined period mapping
$$
\Phi : H^1_{dR}(S^*) \to \C^N
$$
sending a $1$-form $\omega$ to the vector of integrals over these curves, where $N = 2g + (n-1) + (\sum_i (d_i - 1)) + n$.

\medskip

To prove Theorem \ref{mainthm1}, it suffices to prove the following two propositions:

\medskip

\begin{prop} \label{globalinj} The period mapping $\Phi : H^1_{dR}(S^*) \to \C^N$ is injective.
\end{prop}

\medskip

\begin{prop} \label{globalsurj} The period mapping $\Phi : H^1_{dR}(S^*) \to \C^N$ maps onto the
codimension one subspace
$W = \{ (x_1, \dots, x_N) : x_{N-(n-1)} + x_{N - (n-2)} + \dots + x_N = 0 \}$.
\end{prop}

\medskip

\noindent{\bf Proof of Proposition \ref{globalinj}:} Given $[\omega] \in H^1_{dR}(S^*)$ such that all
periods of $\omega$ vanish, in particular the integrals over $a_1, \dots, a_g, b_1, \dots, b_g$ vanish
and the residues at $p_1, \dots, p_n$ vanish (since the integrals over the curves $\gamma^{(i)}$ vanish),
hence $\omega$ has a primitive $F$ on $S'$ which is meromorphic on $S'$. Moreover since $dF = \omega \in \Omega(S^*)$,
$F$ has a Stolz continuous extension to $S^*$ which we may assume satisfies $F(w^{*(1)}_0) = 0$. Since
the integrals of $\omega$ over the curves $\{ \gamma^{(i)}_j \}$ vanish, for each $i = 1, \dots, n$, by
Theorem \ref{localperiod} there
is a germ $f_i \in \OO_{p_i}([h_i])$ and a constant $c_i$ such that $F = f_i + c_i$ in a punctured neighbourhood
of $p_i$. Since the integrals of $\omega$ over the curves $\beta_i$ vanish, $F(w^{*(i)}_0) = F(w^{*(1)}_0) = 0$ for
all $i$, hence $c_i = F(w^{*(i)}_0) = 0$ for all $i$. Thus $F \in \MM(S^*)$ and $[\omega] = 0$. $\diamond$

\medskip

For the proof of Proposition \ref{globalsurj}, we will make use of the following Mergelyan
type Theorem for compact Riemann surfaces due to Gusman (\cite{gusman1}):

\medskip

\begin{theorem} \label{gusman} Let $S$ be a compact Riemann surface and $E \subset S$ a closed
subset such that $S - E$ has finitely many connected components $V_1, \dots, V_m$,
and for each $i$ let $q_i$ be a point of $V_i$. Then any continuous function $f$ on $E$
which is holomorphic in the interior of $E$ can be uniformly approximated on $E$ by functions
meromorphic on $S$ with poles only in the set $\{q_1, \dots, q_m\}$
\end{theorem}

\medskip


We fix an auxiliary point $q$ in the interior
of the polygon $P$ and a small circle $\gamma$ around $q$ such that $q, \gamma$ are disjoint from all
curves in $\CC$. We let $\CC'_4 = \CC_4 \cup \{\gamma\}, \CC' = \CC \cup \{\gamma\}$, and let $\Omega_{II,q}(S^*)$
be the space of all $1$-forms in $\Omega(S^*)$ whose residues are zero at all points of $S'' := S - \{p_1, \dots, p_n, q\}$.
We define another
period mapping
$$
\Psi : \Omega_{II, q}(S^*) \to \C^{N+1}
$$
taking a $1$-form $\omega$ to the vector of its integrals over all curves in $\CC'$.

\medskip

From Theorem \ref{existence}, we can choose a function $f_0 \in \OO(S^*) - \{0\}$. Multiplying $f_0$
by a meromorphic function $g$ which is holomorphic on $S'$ if necessary, we can assume that near each $p_i$ the
function $f_0$ is of the form $f_0 = g_i e^{h_i}$
where $g_i$ is meromorphic near $p_i$ and has a pole at $p_i$ (existence of such a meromorphic function $g$ holomorphic
on $S'$ with poles of high enough order at the points $p_i$ follows from the Riemann-Roch Theorem). We can also
assume that the (finitely many) zeroes and critical points of $f_0$ on $S'$ are disjoint from all curves in $\CC'$
and from $\{q\}$. Then the differential
$\omega_0 := df_0$ lies in $\Omega^{0}(S^*)$, is holomorphic and non-zero on all curves in $\CC'$ and at $q$, and
near each $p_i$ takes the form $\omega_0 = k_i e^{h_i} dw$ where $w$ is a local coordinate near $p_i$ and
$k_i = dg_i/dw + g_i dh_i/dw$ is a meromorphic function with a pole of order at least $d_i + 2$ at $p_i$.

\medskip

We fix a compact subset $E$ given by the union of the curves in $\CC_1, \CC_2$ and
small pairwise disjoint closed disks $U_i$ around the points $p_i$ such that $U_i$ contains the curves
$\{\gamma^{(i)}_j : 1 \leq j \leq d_i - 1 \}$ and $\gamma^{(i)}$. Then $S - E$ is connected and contains the
point $q$.

\medskip

\begin{lemma} \label{typeI} Given $c \in \CC_1$ and $\epsilon > 0$, there exists
$\omega \in \Omega_{II,q}(S^*)$ such that $|\int_c \omega - 1| < \epsilon$
and $|\int_{c'} \omega| < \epsilon$ for all $c' \in \CC' - \{c\}$.
\end{lemma}

\medskip

\noindent{\bf Proof:} Let $f$ be a continuous function on $E$ such that $f$ is supported in a small
arc in the interior of the curve $c \in \CC_1$, $\int_c f\omega_0 = 1$, and $f$ is zero elsewhere on $E$.
Then by Theorem \ref{gusman}, we can find a meromorphic function $R$ on $S$
holomorphic on $S - \{q\}$
such that the $1$-form $\omega := R \omega_0 \in \Omega_{II,q}(S^*)$ satisfies
$|\int_c \omega - 1| < \epsilon$ and $|\int_{c'} \omega| < \epsilon/n$ for all $c' \in \CC - \{c\}$. It
follows that $|\int_{\gamma} \omega| < \epsilon$ since the sum of residues of $\omega$ at $p_1, \dots, p_n, q$
is zero. $\diamond$

\medskip

\begin{lemma} \label{typeII} Given $c \in \CC_2$ and $\epsilon > 0$, there exists
$\omega \in \Omega_{II,q}(S^*)$ such that $|\int_c \omega - 1| < \epsilon$
and $|\int_{c'} \omega| < \epsilon$ for all $c' \in \CC' - \{c\}$.
\end{lemma}

\medskip

\noindent{\bf Proof:} Let $f$ be a continuous function on $E$ which is supported on a small
arc of $c$ disjoint from the disks $U_1, \dots, U_n$, such that $\int_c f\omega_0 = 1$.
Again by Theorem \ref{gusman}, we can find a meromorphic function $R$ on $S$
holomorphic on $S - \{q\}$ such that the $1$-form
$\omega := R \omega_0 \in \Omega_{II,q}(S^*)$ satisfies
$|\int_c \omega - 1| < \epsilon$ and $|\int_{c'} \omega| < \epsilon/n$ for all $c' \in \CC - \{c\}$. It
follows that $|\int_{\gamma} \omega| < \epsilon$ since the sum of residues of $\omega$ at $p_1, \dots, p_n, q$
is zero. $\diamond$

\medskip

\begin{lemma} \label{typeIII} Given $c \in \CC_3$ and $\epsilon > 0$, there exists
$\omega \in \Omega_{II,q}(S^*)$ such that $|\int_c \omega - 1| < \epsilon$
and $|\int_{c'} \omega| < \epsilon$ for all $c' \in \CC' - \{c\}$.
\end{lemma}

\medskip

\noindent{\bf Proof:} Let $c = \gamma^{(i)}_j$ for some $1 \leq i \leq n, 1 \leq j \leq d_i - 1$. By
Proposition \ref{localbasis} we can choose a $1$-form $\omega_i \in \Omega_{p_i}([h_i])$ such that
$\omega_i = \phi_i e^{h_i} dw$ for some meromorphic function $\phi_i$ on $U_i$ with only pole on $U_i$
at $p_i$ of order between $1$ and $d_i$ such that $\int_{\gamma^{(i)}_{j}} \omega_i = 1,
\int_{\gamma^{(i)}_{j'}} \omega_i = \int_{\gamma^{(i)}} \omega_i = 0$ for all $j' \neq j$. We then take
$f$ to be a continuous function on $E$ such that $f \omega_0 = \omega_i$ on $U_i$, $f = 0$ on all other disks $U_{i'},
i' \neq i$, $f = 0$ on all curves in $\CC_1$, and such that the integral of $f \omega_0$ over all
curves in $\CC_2$ is zero. Note that $f$ is holomorphic on the disc $U_i$ since $\phi_i$ has a pole of
order at most $d_i$ at $p_i$ and $k_i$ has a pole of order at least $d_i + 1$ at $p_i$.

Again by Theorem \ref{gusman}, we can find a meromorphic function $R$ on $S$
holomorphic on $S - \{q\}$
such that the $1$-form $\omega := R \omega_0 \in \Omega_{II,q}(S^*)$ satisfies
$|\int_c \omega - 1| < \epsilon$ and $|\int_{c'} \omega| < \epsilon/n$ for all $c' \in \CC - \{c\}$. It
follows that $|\int_{\gamma} \omega| < \epsilon$ since the sum of residues of $\omega$ at $p_1, \dots, p_n, q$
is zero. $\diamond$

\medskip

\begin{lemma} \label{typeIV} Given $y_1, \dots, y_n \in \C$ such that $y_1 + \dots + y_n = 0$
and $\epsilon > 0$, there exists
$\omega \in \Omega_{II,q}(S^*)$ such that $|\int_{\gamma^{(i)}} \omega - y_i| < \epsilon$
for $i = 1, \dots, n, |\int_{\gamma} \omega| < \epsilon$ and $|\int_{c} \omega| < \epsilon$
for all other curves in $\CC'$.
\end{lemma}

\medskip

\noindent{\bf Proof:} As in the proof of the previous Lemma, for each $i$ we choose a $1$-form
$\omega_i \phi_i e^{h_i} dw_i \in \Omega_{p_i}([h_i])$ (where $w_i$ is a local coordinate at $p_i$)
such that the integrals of $\omega_i$ over the curves $\gamma^{(i)}_j$ vanish while the integral
of $\omega_i$ over the circle $\gamma^{(i)}$  equals $y_i$. Let $f$ be a continuous function on $E$ such that
$f \omega_0 = \omega_i$ on $U_i$ (then as before $f$ is holomorphic on $U_i$), $f = 0$ on all
curves in $\CC_1$, and the integrals of $f$ over all curves in $\CC_2$ vanish.

Again by Theorem \ref{gusman}, we can find a meromorphic function $R$ on $S$
holomorphic on $S - \{q\}$
such that the $1$-form $\omega := R \omega_0 \in \Omega_{II,q}(S^*)$ satisfies
$|\int_{\gamma^{(i)}} \omega - y_i| < \epsilon/n$ and $|\int_{c} \omega| < \epsilon/n$ for all other curves
$c \in \CC$. It
follows that $|\int_{\gamma} \omega| < \epsilon$ since the sum of residues of $\omega$ at $p_1, \dots, p_n, q$
is zero. $\diamond$

\medskip

\noindent{\bf Proof of Proposition \ref{globalsurj}}: It follows from Lemmas \ref{typeI}-\ref{typeIV} that
the image of the period mapping $\Psi$ is dense in the codimension two subspace $W'$ of $\C^{N+1}$ given by
$W' := \{ (x_1, \dots, x_{N+1}) \in \C^{N+1} | x_{(N+1) - (n-1)} + \dots + x_{(N+1) - 1} = 0, x_{N+1} = 0 \}$,
hence $W'$ is contained in the image of $\Psi$. Let $V := \Psi^{-1}(W') \subset \Omega_{II, q}(S^*)$, then it
follows from the second equation defining $W'$ that the residue of any $\omega$ in $V$ at the point $q$ vanishes,
thus $V \subset \Omega_{II}(S^*)$ and $\Psi(V) = W'$. It follows that $\Phi(V) = W \subset \C^{N}$. $\diamond$.

\medskip

Theorem \ref{mainthm1} follows immediately from Propositions \ref{globalinj} and \ref{globalsurj}.

\medskip

For the proof of Theorem \ref{mainthm2}, we consider a similar period mapping (keeping the
same notation as before)
$$
\Phi : H^1_{dR, 0}(S^*) \to \C^N.
$$
As before it suffices to show that $\Phi$ is injective and that $\Phi$ surjects onto the
codimension one subspace $W$. The proof of injectivity is the same as before, observing that if
$F \in \MM(S^*)$ is the primitive of a $1$-form $\omega \in \Omega^{0}(S^*)$ then
$F \in \OO(S^*)$. We proceed to the proof of surjectivity.

\medskip

As before we fix the point $q$
and the function $f_0$ and $1$-form $\omega_0 = df_0$. The following lemma is a
straightforward consequence of the Riemann-Roch Theorem:

\medskip

\begin{lemma} \label{Rpq} Given an integer $l$ and a point $p \neq q$,
there exists a meromorphic function $R$ on $S$ such that $R$ is holomorphic
on $S - \{p,q\}$ and the order of $R$ at $q$ is equal to $-l$.
\end{lemma}

\medskip

\noindent{\bf Proof:} For $N$ a large enough positive integer, there exists $R \in \OO_{D} - \OO_{D'}$,
where $D = p^{N}q^{l+1}, D' = p^{N}q^{l}$. $\diamond$

\medskip

\begin{lemma} \label{addexact} Given $\omega \in \Omega(S^*)$ such that $\omega$ is holomorphic on $S'' = S - \{p_1, \dots, p_n, q\}$,
there exists $f \in \MM(S^*)$ such that $f$ is holomorphic on $S''$ and $\omega + df$ has at most a simple
pole at $q$.
\end{lemma}

\medskip

\noindent{\bf Proof:} Let $\alpha$ be the logarithmic differential $df_0/f_0$,
which is a meromorphic $1$-form on $S$, and fix a local coordinate $z$ at the point $q$
such that $z(q) = 0$. We may assume that $\omega$ has a pole at $q$ of order $m \geq 2$
(otherwise we can take $f = 0$ and we are done).
Then, $f_0$ is non-zero at $q$, so near $q$ we have
$$
\omega = \left( \frac{b_{-m}}{z^m} + \dots + \frac{b_{-1}}{z} + O(1) \right) \cdot f_0 \cdot dz
$$
for some constants $\{b_k\}$. By Lemma \ref{Rpq} we can choose a function $R$ meromorphic on $S - \{p_1, q\}$
such that $R$ has a pole of order $(m-1)$ at $q$, and multiplying by a constant if necessary we may assume
that
$$
R(z) = \frac{c_{-(m-1)}}{z^{m-1}} + O\left(\frac{1}{z^{m-2}}\right)
$$
near $q$, where $c_{-(m-1)} = -b_{-m}/(m-1)$. Then $R \cdot f_0 \in \MM(S^*)$ is holomorphic on $S''$,
and $d(R f_0) = (dR + R \alpha)f_0$, where $\alpha$ is holomorphic near $q$, so
\begin{align*}
\omega + d(Rf_0) & = \left[ \left( \frac{b_{-m}}{z^m} + O\left(\frac{1}{z^{m-1}}\right) \right) +
\left( \frac{-(m-1)c_{-(m-1)}}{z^m} + O\left(\frac{1}{z^{m-1}}\right) \right) + O\left(\frac{1}{z^{m-1}}\right) \cdot O\left(1\right)
\right] \cdot f_0 \cdot dz \\
                 & = O\left(\frac{1}{z^{m-1}}\right) \cdot f_0 \cdot dz \\
\end{align*}

Thus $\omega + d(Rf_0)$ has at most a pole of order $(m-1)$ at $q$, so we are done by induction on $m$.
$\diamond$

\medskip

\noindent{\bf Proof of \ref{mainthm2}:} Let $\Psi : \Omega_{II, q}(S^*) \to \C^{N+1}$ be as
before.

\medskip

It follows from the proof of Proposition \ref{globalsurj},
that for any vector $\vec{v} \in W'$ and any $\epsilon > 0$ there is an $\omega \in \Omega(S^*)$
which is holomorphic on $S''$ such that $||\Psi(\omega) - \vec{v}|| < \epsilon$. By the previous lemma,
we can find a $1$-form $\omega' = \omega + df$ such that $\omega' \in \Omega(S^*)$, $\omega'$ is holomorphic
on $S''$, has at most a simple pole at $q$, and $\Psi(\omega') = \Psi(\omega)$. Letting $V_q$ denote
the subspace of $\Omega(S^*)$ consisting of $1$-forms which are holomorphic on $S''$ with at most a
simple pole at $q$, it follows that $\Psi(V_q)$ is dense in, and hence contains, $W'$. Letting $V' \subset V_q$ be the subspace
of $V_q$ such that $\Psi(V') = W'$, we have $V' \subset \Omega^{0}(S^*)$ (since $1$-forms in $V'$ have
at most a simple pole at $q$ and their residue at $q$ vanishes) and $\Phi(V') = W$, so
$\Phi : H^1_{dR, 0}(S^*) \to W$ is surjective. $\diamond$

\bigskip

%
%
%
%
%
%
%
%

\bibliography{expsingular}
\bibliographystyle{alpha}

\end{document}